\documentclass[11pt]{amsart}
\usepackage{amssymb,amsfonts,amsxtra}
\usepackage[all]{xy}
\usepackage{fullpage}
\DeclareMathOperator{\Der}{Der}
\DeclareMathOperator{\Coder}{Coder}
\DeclareMathOperator{\Hom}{Hom}

\DeclareMathOperator{\im}{Im}

\def\lla{\longleftarrow}

\def\lra{\longrightarrow}

\newtheorem{theorem}{Theorem}[section]
\newtheorem{cor}[theorem]{Corollary}
\newtheorem{lem}[theorem]{Lemma}
\newtheorem{prop}[theorem]{Proposition}

\theoremstyle{definition}
\newtheorem{defi}[theorem]{Definition}
\newtheorem{rem}[theorem]{Remark}
\numberwithin{equation}{section} 
\thanks{J.B. was partially supported by the NSF grant DMS-0204558}
\thanks{A.L. was partially supported by the EPSRC grant No. GR/R84276/01}
\begin{document}
\title[Andr\'e-Quillen cohomology and rational homotopy of function spaces]
{Andr\'e-Quillen cohomology and rational homotopy of function
spaces}
\author{J.Block \& A. Lazarev}
\address {Department of Mathematics, Bristol University, Bristol, BS8 1TW,
England.} \email{a.lazarev@bristol.ac.uk}
\address{Department of Mathematics,
University of Pennsylvania, Philadelphia, PA 19104, USA.}
\email{blockj@math.upenn.edu}
\begin{abstract}
We develop a simple theory of Andr\'e-Quillen cohomology for
commutative differential graded algebras over a field of
characteristic zero. We then relate it to the homotopy groups of
function spaces and spaces of homotopy self-equivalences of
rational nilpotent $CW$-complexes. This puts certain results of
Sullivan in a more conceptual framework.
\end{abstract} \maketitle
\section{Introduction}
Andr\'e-Quillen cohomology is a cohomology theory for commutative
algebras originally introduced in \cite{Andre} and \cite{Qui}. It
was subsequently generalized to cover simplicial algebras over
operads, \cite{GH}, differential graded $E_\infty$-algebras,
\cite{Man} and commutative $S$-algebras, \cite{Bas}.

One of the purposes of the present paper is to give a simple and
direct treatment of the Andr\'e-Quillen cohomology in the category
of commutative differential graded algebras (dga's for short) over
a field of characteristic zero. This is done in Section 2. Our
initial definition of Andr\'e-Quillen cohomology of a dga $A$ with
coefficients in the differential graded (dg) module $M$ over $A$
is via an explicit cochain complex $C^*_{AQ}(A,M)$ similar to the
one introduced in \cite{Har}. We then produce various equivalent
characterizations of Andr\'e-Quillen cohomology, introduce the
Gerstenhaber bracket on $C^*_{AQ}(A,A)$ and show its homotopy
invariance. In this connection we mention the recent paper
\cite{FMT} where analogues of some of our results were proved in
the context of Hochschild cohomology.

In Section 3 we apply the developed techniques to computing the
homotopy groups of function spaces. (We are dealing with
\emph{unpointed} spaces, however our machinery could be easily
adapted to the pointed situation as well.) In particular we are
concerned with the group $hAut(X)$ of homotopy classes of homotopy
self-equivalences of a nilpotent $CW$-complex $X$. A well-known
theorem of Sullivan \cite{Sul} and Wilkerson \cite{Wil} asserts
that under suitable finiteness assumptions $hAut(X)$ is an
arithmetic group, that is, commensurable to the group of integer
points of some algebraic group over $\mathbb{Q}$. An important
step is to show that the group $hAut(X_{\mathbb{Q}})$ is
isomorphic to the group of $\mathbb{Q}$-points of an algebraic
group. Here $X_{\mathbb{Q}}$ denotes the \emph{rationalization} of
the space $X$, i.e. its localization with respect to the homology
theory $H_*(-,\mathbb{Q})$.

We reprove this result and identify the Lie algebra of this
algebraic group. It turns out to be isomorphic to
$H^0_{AQ}(A^*(X),A^*(X))$, the zeroth Andr\'e-Quillen cohomology
of the Sullivan-de Rham algebra of $X$ with coefficients in
itself. The Lie bracket corresponds to the Gerstenhaber bracket on
$H^*_{AQ}(A^*(X),A^*(X))$.

We also consider the question of computing the higher homotopy
groups of a function space $F(X,Y)$ for two rational spaces $X$
and $Y$. The answer is again formulated in terms of
Andr\'e-Quillen cohomology associated to the Sullivan-de Rham
models of $X$ and $Y$. This result was hinted at in \cite{Laz}.

\section{Andr\'e-Quillen cohomology of commutative differential graded algebras.}
Let $\mathcal C$ denote the category of commutative differential
graded algebras over $k$, not necessarily connected. Here $k$ is a
field of characteristic $0$. The differential is assumed to
\emph{raise} the degree. Then $\mathcal C$ admits a structure a a
closed model category as follows:
\begin{itemize}\item weak equivalences are maps which induce an
isomorphism on cohomologies; \item fibrations are surjective
maps;\item cofibrations are the maps which have the left lifting
property with respect to acyclic fibrations.\end{itemize} That
$\mathcal C$ is a closed model category is proved in the case of
connected dga's in \cite{BG}. The general case is due to Hinich,
\cite{H} who proved it in the still greater generality of algebras
over an operad.

Let us describe the cofibrant objects in $\mathcal C$. First
consider the operation of glueing cells to a dga (also called the
Hirsch extension). Let $A$ be a dga and $V$ be a graded vector
space. Let $f:V\lra Z^*(A)$ be a linear map of degree $1$ from $V$
to the space of cocycles of $A$. Then define a new dga $A_f$ whose
underlying graded vector space is $A\otimes \Lambda V$. Here we
denoted, following tradition, by $\Lambda V$ the free graded
commutative algebra on the vector space $V$. The differential on
$A$ is the old differential and the one on $V$ is given by the map
$f$. Then $A_f$ is said to be obtained from $A$ by glueing a
(generalized) cell. Observe that $A_f$ is a pushout of $A$ by a
free commutative algebra which justifies the name. A dga obtained
from the trivial dga $k$ is called a \emph{cell} dga. Any cell dga
is cofibrant and any cofibrant dga is a retract of a cell dga.

Now consider the category $A-mod$ of dg modules over a dga $A$.
This is also a closed model category where fibrations are
surjective maps.  Then a graded derivation of $A$ with values in
$M$ of degree $d$ is a map $\xi:A\lra M$ of degree $d$ which
satisfies the Leibnitz condition:
\[\xi(b_1b_2)=\xi(b_1)b_2+(-1)^{d|b_1|}b_1\xi(b_2).\] The collection of
all derivations form a complex, in fact a dg $A$-module
$\Der^*(A,M)$.

Associated to a dga $A$ is the dg module of its K\"ahler
differentials $\Omega_A$. It is defined in the usual manner as
$\Omega_A:=I/I^2$ where $I$ is the kernel of the multiplication
map $A\otimes A\lra A$. It is a standard fact that there is an
isomorphism of dg $A$-modules:
\[\Hom^*_A(\Omega_A,M)\cong \Der^*(A,M).\]
Let us now introduce the derived version of $\Omega_A$ also called
the \emph{Andr\'e-Quillen homology} of $A$. First recall that the
(homological) Hochschild complex of $A$ with coefficients in
itself is defined as the complex
\[C_*(A,A)=\{A\lla  A^{\otimes 2}\lra\ldots\}\]
with the standard bar differential. $C_*(A,A)$ is in fact itself a
differential graded algebra with respect to the shuffle product,
since $A$ is graded commutative. Since $A$ is a dga this is in
fact a \emph{bicomplex}. (The Hochschild differential lowers
degree while the differential of $A$ raises degree. Thus the total
degree in the bicomplex is the difference of the two.) We will
make use of the truncated version of $C_*(A,A)$ denoted by
$\bar{C}_*(A,A)$. This is the same complex as $C_*(A,A)$ but
starting with $A^{\otimes 2}$. Since $A$ is commutative the
complex $\bar{C}_*(A,A)$ splits off ${C}_*(A,A)$ as a direct
summand.

\begin{defi}The Andr\'e-Quillen complex $C^{AQ}_*(A,A)$ of a dga $A$ is
the quotient complex of $\bar{C}_*(A,A)$ by the subcomplex of
decomposables, i.e. those elements which could be represented as
shuffle products of two or more elements in $\bar{C}_*(A,A)$.
\end{defi}
\begin{rem} This complex (shifted) was defined by Harrison in \cite{Har} in the
case when $A$ is a usual (ungraded) algebra. Its homology is also
called the Harrison homology of $A$. It is well known that in
characteristic zero case the shifted Harrison homology agree with
the Andr\'e-Quillen homology defined by means of simplicial
resolutions.
\end{rem}

\begin{theorem}
Let $A$ be a cofibrant dga. Then there is a quasi-isomorphism of
dg $A$-modules:
\[\Omega_A\simeq C^{AQ}_*(A,A).\]\end{theorem}
\begin{proof} Consider the map $f:A^{\otimes 3}\lra A^{\otimes 2}:$
\[f:a\otimes b\otimes c\mapsto ab\otimes c-a\otimes
bc+(-1)^{(|a|+|b|)|c|}ca\otimes
 b.\] Clearly, $\im f=\Omega_A$.
There results a map of dg modules
\[C^{AQ}_*(A,A)\lra \Omega_A\] and we want to prove that
this map is a quasi-isomorphism for a cofibrant commutative dga
$A$.

Without loss of generality we assume that $A$ is constructed from
$k$ by a series of Hirsch extensions. This gives $A$ a filtration
for which the associated graded algebra is simply the free
commutative algebra on some set of generators with zero
differential. This filtration lifts to $\Omega_A$ and
$C^{AQ}_*(A,A)$ so that the canonical map $C^{AQ}_*(A,A)\lra
\Omega_A$ is a filtered map. Since it is clearly a
quasi-isomorphism if $A$ is free commutative with vanishing
differential we conclude that the map is a quasi-isomorphism on
the level of associated graded modules, therefore it was a
quasi-isomorphism to begin with.

\end{proof}

Now let us turn to the functor of derivations and its derived
version. Let $M$ be a dg module over a dga $A$ and denote by
$\tilde{A}$ the cofibrant replacement of $A$. Then $M$ is also a
dg $\tilde{A}$-module and we define the derived functor of
derivations of $A$ with values on $M$ as
 $\Der^*(\tilde{A},M)$. (Strictly speaking, we have not set
up things so that it is a derived functor; Der is not even a
functor.) We have \[\Der^*(\tilde{A},M)\cong
\Hom^*_A(\Omega_{\tilde{A}},M)\simeq \Hom^*_A(C^{AQ}_*(A,A),M).\]
The complex $\Hom^*_A(C^{AQ}_*(A,A),M)$ is the subcomplex of the
truncated Hochschild complex
\[\bar{C}^*(A,M):=\Hom^*_A(\bar{C}_*(A,A),M)\] consisting of those
cochains which vanish on the shuffle products. (A shifted version
of this complex is commonly called the Harrison cohomology complex
of $A$ with coefficients in $M$.) We denote this complex by
$C^*_{AQ}(A,M)$ and its cohomology by $H^*_{AQ}(A,M)$. Therefore
we proved the following
\begin{theorem}\label{derived}
The cohomology of the complex $C^*_{AQ}(A,M)$ is isomorphic to the
cohomology of the differential graded module $\Der^*(\tilde{A},M)$
where $\tilde{A}$ is a cofibrant replacement of $A$.
\end{theorem}
\begin{cor}Let $A$ be a cofibrant dga and $M$ is a dg $A$-module. Then there is a
spectral sequence $H^*_{AQ}(H^*(A),H^*(M))\Longrightarrow
H^*_{AQ}(A,M)$.\end{cor}

Thus, we have two ways to compute the Andr\'e-Quillen cohomology
of a dga $A$ with values in a dg $A$-module $M$. The first is to
replace $A$ with its cofibrant approximation and take its
derivations in $M$. The second is via the functorial complex
$C^*_{AQ}(A,M)$. The method via the complex $C^*_{AQ}(A,M)$ is
better suited for theoretical purposes; in particular it gives
rise to the spectral sequence as above. Another useful property of
the complex $C^*_{AQ}(A,M)$ is that it is a direct summand of the
Hochschild complex. This is something that is not seen from the
point of view of the derived functor of the derivations.

Derivations also admit the following useful interpretation in
terms of square-zero extensions. Let $A$ be a dga and $M$ be a dg
$A$-module. Denote by $A\ltimes M$ the dga which is isomorphic as
a complex to $A\oplus M$ with multiplication defined as
\[(a_1,m_1)(a_2,m_2)=(a_1a_2,a_1m_2+m_1a_2).\] We will call
$A\ltimes M$ the \emph{square-zero extension} of $A$ by $M$. The
dga $A\ltimes M$ is supplied with a dga map into $A$ which is
simply the projection onto the first component. We can thus view
$A\ltimes M$ as an object in the \emph{overcategory} of $A$, i.e.
the category whose objects are dga's $B$ supplied with a map
$\epsilon_B:B\lra A$ and morphisms are obvious commutative
triangles. We will denote this overcategory by $\mathcal{C}_A$. It
inherits the structure of a closed model category from
$\mathcal{C}$ so that a morphism in $\mathcal{C}_A$ is a
cofibration if it is so considered as a map in $\mathcal{C}$. Note
that $A$ is also an object in $\mathcal{C}_A$ in an obvious
fashion.

The association $M\rightsquigarrow A\ltimes M$ is a functor
$A-mod\mapsto \mathcal{C}_A$. It clearly preserves weak
equivalences and therefore lifts to a functor between the
corresponding homotopy categories.

Let $B$ be an object in $\mathcal{C}_A$ which we could assume to
be cofibrant without loss of generality. Then an $A$-module $M$
has a structure of a $B$-module via the structure map
$\epsilon_B:B\lra A$. An elementary calculation shows that a
derivation $B\lra M$ is nothing but a map $B\lra A\ltimes M$ in
$\mathcal{C}_A$. Furthermore, we have the following proposition:
\begin{prop}\label{characterization}There is a natural isomorphism
\[H^0_{AQ}(B,M)\cong [B,A\ltimes M]_{\mathcal{C}_A}\]
where $[-,-]_{\mathcal{C}_A}$ denotes the set of homotopy classes
of maps in $\mathcal{C}_A$.\end{prop}
\begin{proof}The map $b\lra b\otimes 1-1\otimes b:B\lra B\otimes
B$ can be considered as a map $B\lra I$ where $I$ is the kernel of
the multiplication map $B\otimes B\lra B$. Composing this map with
the projection $I\lra I/I^2=\Omega_B$ we get a map $d_B:B\lra
\Omega_B$, known as the \emph{universal derivation} of the dga
$B$. (It is standard to check that this is indeed a derivation).
Denoting by $[-,-]_{B-mod}$ the set of morphisms in the homotopy
category of $B-mod$ we get a natural transformation
\begin{equation}\label{AQ^0}H_{AQ}^0(B,M)=[\Omega_B,M]_{B-mod}\mapsto [B,A\ltimes
M]_{\mathcal{C}_A}\end{equation} which associates to a homotopy
class of a map $f:\Omega_B\lra M$ the composition
\[\xymatrix{B\ar^{(1,d_B)}[r] &B\ltimes
\Omega_B\ar^{(\epsilon_B,f)}[r]&A\ltimes M.}\] To check that the
map of (\ref{AQ^0}) is an isomorphism consider the functor
$\mathbb{S}$ which associates to a complex $N$ the symmetric
algebra $\mathbb{S}(N)$. The functor $\mathbb{S}$ is left adjoint
to the forgetful functor from dga's to complexes and this
adjunction passes to homotopy categories. Now it is easy to see
that the map of \ref{AQ^0} is an isomorphism if $B=\mathbb{S}(N)$.
The general case follows by virtue of the following canonical
split coequalizer which exists for any dga $B$:
\[\xymatrix{\mathbb{S}^2B\ar@<-0.5ex>[r]\ar@<0.5ex>[r]&\mathbb{S}B\ar[r]&B.}\]
\end{proof}
\begin{rem}On a technical note observe that the canonical
projection $A\ltimes M\lra A$ is a fibration from which it follows
that $A\ltimes M$ is a fibrant object in $\mathcal{C}_A$.
Therefore the set $[B,A\ltimes M]_{\mathcal{C}_A}$ does represent
the set of morphisms in the homotopy category of
$\mathcal{C}_A$.\end{rem}

The next thing we are going to describe is the \emph{Gerstenhaber
bracket} on the Andr\'e-Quillen cohomology. It turns out that for
a dga $A$ the complex $C^*_{AQ}(A,A)$ admits the structure of a dg
Lie algebra. The most conceptual way to describe it is due to
Schlessinger and Stasheff, cf.\cite{SS} which we will now recall.

Let $CA$ be the free Lie coalgebra on $A$. It could be described
as the quotient of the reduced tensor algebra
$TA_+=\oplus_{i=1}^\infty A^{\otimes i}$ by the image of the
shuffle product. The usual bar differential descends from $TA_+$
to $CA$ making it an acyclic Lie coalgebra. Then $C^*_{AQ}(A,A)$
is naturally identified with the space $\Coder^*(CA,CA)$ of
coderivations of the Lie coalgebra $CA$. Moreover,
$\Coder^*(CA,CA)$ is naturally a dg Lie algebra with respect to
the commutator bracket which we will call the Gerstenhaber bracket
since it is a direct analogue of the bracket introduced in
\cite{Ger} on the Hochschild complex of an algebra.

 On the other hand Andr\'e-Quillen
cohomology of $A$ could be described as
$\Der^*(\tilde{A},\tilde{A})$ where $\tilde{A}$ is the cofibrant
approximation of $A$. This gives another structure of a graded Lie
algebra on $H_{AQ}^*(A,A)$. The following result shows that these
two structures coincide, and, moreover are invariants of the weak
homotopy type of $A$.
\begin{theorem}\begin{enumerate}\item
For two weakly equivalent cofibrant dga's $A$ and $B$ the dg Lie
algebras $\Der^*(A,A)$ and $\Der^*(B,B)$ are quasi-isomorphic;
\item for two weakly equivalent (not necessarily cofibrant) dga's
$A$ and $B$ the dg Lie algebras $\Coder^*(CA,CA)$ and
$\Coder^*(CB,CB)$ are quasi-isomorphic;\item for a cofibrant dga
$A$ the dg Lie algebras $\Der^*(A,A)$ and $\Coder^*(CA,CA)$ are
quasi-isomorphic.\end{enumerate}\end{theorem}
\begin{proof}The main problem is, of course, that $\Der^*(A,A)$ is not functorial
with respect to $A$. To overcome this difficulty we use the
properties of closed model categories. For part (1) let $f:A\lra
B$ be a weak equivalence between two cofibrant dga's. First assume
that $f$ is a fibration. Then there exists a right splitting
$g:B\lra A$ so that $f\circ g=id_B$. Then define the map of dga's:
$End(B)\lra End(A)$ by assigning to a map $s:B\lra B$ the map
$g\circ s\circ f$. The condition $f\circ g=id_B$ ensures that this
map respects composition of endomorphisms. It follows that we have
a map of dg Lie algebras $h:\Der^*(B,B)\lra \Der^*(A,A)$ together
with the commutative diagram of complexes:
\[\xymatrix{\Der(A,A)\ar[dr]\ar^h[rr]&&\Der(B,B)\\&\Der(A,B)\ar[ur]}\]
Here the southeast arrow and the northeast arrows are induced by
$f$ and $g$ respectively. It follows that $h$ is a
quasi-isomorphism.

 Similarly if $f:A\lra B$ is a
cofibration then it admits a left splitting $g:B\lra A$ so that
$g\circ f=id_A$ and we have a map of dg Lie algebras
$\Der^*(B,B)\lra \Der^*(A,A)$. In the general case we use the
presentation of $f$ as a composition of a cofibration and a
fibration.

The argument for part (2) is more difficult, although the idea is
the same. The most conceptual way is to introduce a closed model
category structure on dg Lie coalgebras and notice that for a dga
$A$ the dg Lie coalgebra $CA$ is a fibrant-cofibrant object. The
corresponding result for dg coalgebras is due to Hinich, cf.
\cite{Hinich}.

While there is little doubt that this could be done a detailed
proof would require rewriting much of Hinich's paper (with
appropriate modifications). Since we don't need the full strength
of the model category structure we give an alternative proof here
which is of some independent interest.

Let $f:A\lra B$ be a quasi-isomorphism between two dga's and
denote by $\tilde{f}$ the induced map of dg Lie coalgebras $CA\lra
CB$. Note that $\tilde{f}$ is a quasi-isomorphism. Moreover, note
that $CA$ and $CB$ have filtrations inherited from the tensor
coalgebras and $\tilde{f}$ is a \emph{filtered} quasi-isomorphism,
that is it induces a quasi-isomorphism on each graded component.

 We claim that  $\tilde{f}$ could be factored as
\[\xymatrix{CA\ar^{\tilde{g}}[r]&D\ar^{\tilde{h}}[r]&CB}\]
where \begin{itemize} \item both maps $\tilde{g}$ and $\tilde{h}$
are filtered quasi-isomorphisms of dg Lie coalgebras; \item $D$ is
a dg Lie coalgebra that is free as a Lie coalgebra; \item the map
$\tilde{g}:CA\lra D$ admits a filtered left inverse;\item the map
$\tilde{h}:D\lra CB$ admits a filtered right inverse.\end{itemize}

This claim clearly implies the result we need by standard spectral
sequence arguments. We will now begin to prove the required
factorization. It would be convenient for us to work with
\emph{complete free} Lie algebras instead of cofree Lie
coalgebras. Let  $LA$ and $LB$ denote the (graded) $k$-linear dual
to $CA$ and $CB$ with their compact-linear topology. (The reader
is referred to \cite{Kot} for an account on linear compactness.)
Thus, $LA$ and $LB$ are complete free Lie algebras on the
$k$-vector space $A^*$ and $B^*$. Moreover, $LA$ and $LB$ have
differentials making them dg Lie algebras and there is a
continuous map $\tilde{f}^*:LB\lra LA$. Further, $LA$ and $LB$
possess natural filtrations by bracket length and $\tilde{f}^*$ is
a a filtered quasi-isomorphism.

Under these conditions we will construct a dg Lie algebra
$D^\prime$ so that $\tilde{f}^*$ factors as
\[\xymatrix{LB\ar^{\tilde{h}^*}[r]& D^\prime\ar^{\tilde{g}^*}[r]& LA}\] so
that:\begin{itemize} \item $\tilde{h}^*$ and $\tilde{g}^*$ are
continuous filtered quasi-isomorphisms of dg Lie algebras;
 \item $D^\prime$ is a dg Lie algebra
which is free and complete as a Lie algebra; \item the map
$\tilde{h}^*:LB\lra D$ admits a continuous filtered left
inverse;\item the map $\tilde{g}^*:D^\prime\lra LA$ admits a
continuous filtered right inverse.\end{itemize} Taking continuous
duals this is easily seen to be equivalent to the statement of the
claim.

We use an analogue of the mapping cylinder construction in
topology. Let $L\langle t\rangle$ denote the completed free Lie
algebra with two generators $t,dt$ such that $|dt|=|t|+1$. We
introduce a differential in $L\langle t\rangle$ by setting
$d(t)=dt$ and $d(dt)=0$. Then $L\langle t\rangle$ is a filtered dg
Lie algebra which is filtered contractible (which means that the
graded components corresponding to its filtration are
contractible). For a filtered dg Lie algebra $l$ the free Lie
product of $l$ and $L\langle t\rangle$ has an induced filtration
and we denote by $l\langle t\rangle$ its completion with respect
to this filtration. For a finite indexing set $I$ we introduce the
notation $l\langle t_\alpha\rangle_{\alpha\in I}$ to denote the
completion of the free Lie product of $l$ and the collection of
Lie algebras $L\langle t_\alpha\rangle$ with $\alpha\in I$. If $I$
is infinite we denote by $l\langle t_\alpha\rangle_{\alpha\in I}$
the inverse limit of $l\langle t_\alpha\rangle_{\alpha\in J}$
where $J$ ranges through finite subsets of $I$. Note that the
natural inclusion $L\hookrightarrow L\langle t_\alpha\rangle$ is a
filtered quasi-isomorphism having a filtered right inverse (which
sends all $t_\alpha$ to zero).

 Now let
\[D^\prime:=LB\langle t_\alpha\rangle\]
where $\alpha$ runs through the set of homogeneous elements of
$LA$. We then have a factorization of $\tilde{f}^*$ as
\[LB\lra D^\prime\lra LA\]
where the first arrow is the obvious inclusion and the second
arrow is the map which coincides with $\tilde f^*$ on $LB$ and
takes each element $t_\alpha$ in $LB\langle t_\alpha\rangle$ to
$\alpha\in LA$. Note that the latter map is a \emph{surjection}
and a filtered quasi-isomorphism. The following lemma ensures that
it has a (filtered) right inverse.
\begin{lem} Let $s:l\lra g$ be a continuous surjective map of
dg Lie algebras which are completions of free Lie algebras on
linearly compact spaces. Assume that $s$ is a filtered
quasi-isomorphism. Then $s$ admits a continuous filtered right
inverse.
\end{lem}
\begin{proof} Choose a filtered Lie algebra map $i:g\lra l$ for which $i\circ
s=id_l$ (we do not claim that $i$ is a map of dg Lie algebras).
This is possible since $g$ is a complete free Lie algebra. Then we
have an isomorphism of vector spaces $l\cong I\oplus i(g)$ where
$I$ is the kernel of $s$. For $a\in g$ we have
\[d(i(a))=i(d(a))+\xi(a).\] It is straightforward to check that
the map $\xi:g\lra I$ is a cocycle in the complex $\Der^*_f(g,I)$
of filtration-preserving derivations of $g$ with values in $I$
(the $g$-module structure on $I$ is provided by the map $i$).
Furthermore, it is possible to find a right inverse to $s$
compatible with differentials if and only if $\xi$ is a coboundary
in $\Der^*_f(g,I)$. Since $s$ is a filtered quasi-isomorphism, $I$
is filtered contractible. Denote by $F_n(I)$ the $n$th filtration
component of $I$ i.e. $F_n(I)$ consists of those elements in $I$
which have bracket length $\geq n$. Then $F_n(I)/F_{n+1}(I)$ has
zero homology. We have \[\Der^*_f(g, I)\cong
\lim_{\leftarrow}\Der^*(g,F_n(I)).\] However since
$F_n(I)/F_{n+1}(I)$ is contractible the obvious induction shows
that $\Der^*(g,F_n(I))$ is likewise contractible so
$\Der^*_f(g,I)$ is contractible and our lemma is
proved.\end{proof}

Finally for (3) we only need to note that the canonical projection
$\Coder^*(CA,CA)\lra \Der^*(A,A)$ is a map of dg Lie algebras. The
required result then follows from Theorem \ref{derived} and part
(2) of the present theorem which has just been proved.

\end{proof}

\section{Rational homotopy groups of function spaces}
Before discussing the rational homotopy of function spaces we need
to establish the following facts about the relation of spaces of
maps and that of their Postnikov stages. These are presumably
well-known but we have been unable to find a proper reference. Let
$hAut(X)$ denote the group of homotopy classes of homotopy
equivalences of a space $X$. The set of homotopy classes of
(unpointed) maps $X\lra Y$ is denoted by $[X,Y]$, as usual.

\begin{prop}\label{wellknown}
Let $X$ and $Y$ be connected CW-complexes with $\dim X=n$ and $Y$
is nilpotent. Let $X_n$ and $Y_n$ denote their $n$th Postnikov
stages. Then
\begin{enumerate} \item   $[X, Y]\lra [X_n,Y_n]$ is a
bijection. \item If $X$ is nilpotent then $hAut(X)\lra hAut(X_n)$
is an isomorphism of groups.
\end{enumerate}
\end{prop}
\begin{proof}
Note first that part (2) of the Theorem is an immediate
consequence of part (1).

For (1) consider the following maps of sets:
\[[X,Y]\lra[X,Y_n]\lla[X_n,Y_n].\]
Let us prove first that $[X,Y]\lra[X,Y_n]$ is a bijection. For
simplicity we assume that the Postnikov tower of $Y$ consists of
principal fibrations. In the general case we could argue
similarly, replacing the Postnikov tower of $Y$ by its principal
refinement (which exists since $Y$ is nilpotent). A map $X\lra
Y_n$ lifts to a map $X\lra Y_{n+1}$ if and only if an obstruction
class lying in $H^{n+2}(X,\pi_{n+1}(Y))$ is zero. This is ensured
by our assumption that the dimension of $X$ is less than or equal
to $n$. Furthermore, the group $H^{n+1}(X,\pi_{n+1}(Y))$ acts on
the set $[X,Y_{n+1}]$ so that the set of orbits is precisely
$[X,Y_n]$. Since this group is zero we conclude that there is a
bijection between sets $[X,Y_n]$ and $[X,Y_{n+1}]$. Arguing by
induction up the Postnikov tower of $Y$ we see that
$[X,Y]\lra[X,Y_n]$ is a bijection.

We will now show that the map $[X,Y_n]\lla[X_n,Y_n]$ is bijective.
First of all, since associating to a space its $n$th Postnikov
stage is a functor in the homotopy category it follows that any
map $X\lra Y_n$ extends to $X_n\lra Y_n$. Furthermore, assuming
again that the Postnikov tower of $Y$ consists of principal
fibrations we see that the ambiguities in choosing extensions lie
in the relative cohomology groups $H^k(X,X^n,\pi_k(Y_n))$. For
$k\leq n$ these groups vanish since $X$ is an $n$-dimensional
complex whereas for $k>n$ they vanish since $Y_n$ has no homotopy
above dimension $n$. Therefore any map $X\lra Y_n$ extends to
$X_n\lra Y_n$ uniquely up to homotopy.
\end{proof}
We will also need the following linearized version of the
preceding result.
\begin{prop}\label{wellknowntoo}Let $A$ and $B$ be two connected dga's so that the
cohomology  of $A$ vanishes above dimension $n$ and $A\lra B$ be a
fixed dga map making $B$ into an $A$ dg-module. Let $A_n$ and
$B_n$ denote the $n$th Postnikov stages of $A$ and $B$
respectively. Then
\begin{enumerate}\item $H_{AQ}^0(A, B)\lra H_{AQ}^0(A_n,B_n)$ is a
bijection. \item $H_{AQ}^0(A, B)\lra H_{AQ}^0(A_n,B_n)$ is an
isomorphism of Lie algebras.
\end{enumerate}
\end{prop}
\begin{proof} Let us first explain the construction of the map
$H_{AQ}^0(A, B)\lra H_{AQ}^0(A_n,B_n)$ figuring in the formulation
of the theorem. We assume from the very beginning that $A$ and $B$
are minimal which results in no loss of generality. Then
$H_{AQ}^0(A,B)$ is simply the zeroth cohomology of the complex
$\Der^*(A,B)$ (this, of course, is not a dg Lie algebra unless
$A=B$). Next, $A_n$ and $B_n$ are subalgebras of $A$ and $B$
generated by the polynomial generators in degrees $\leq n$.
Because of the minimality $A_n$ and $B_n$ are closed under the
differential. Clearly a derivation of degree $0$ $A\lra B$
descends to a derivation $A_n\lra B_n$. Moreover, cycles in
$\Der^0(A,B)$ map to cycles in $\Der^0(A_n,B_n)$ and boundaries -
to boundaries. Thus, the map $H_{AQ}^0(A, B)\lra
H_{AQ}^0(A_n,B_n)$ is well-defined and clearly is a Lie algebra
map in the case $A=B$. It is further obvious that part (1) of the
proposition is a consequence of part (2).

As in the proof of Proposition \ref{wellknown} we assume, purely
for notational simplicity, that the Postnikov tower of $A$
consists of principal fibrations. That means that for a polynomial
generator $x$ of $A$ of degree $k$ the element $dx$ belongs to
$A_{k-1}$. If this is not the case, then one could argue
similarly, using the principal refinement of the Postnikov tower
of $A$ which exists since $A$ is minimal.

Let us denote by $d_A$ and $d_B$ the differentials in $A$ and $B$
and, by abuse of notation, also the differentials in $A_n$ and
$B_n$.

Surjectivity: let $\xi:A_n\lra B_n$ be a derivation of zero degree
such that $\xi\circ d_B=d_A\circ \xi$. In other words, $\xi$ is a
cycle in $\Der^0(A_n,B_n)$. Then $\xi$ could be extended to
$A_{n+1}$ if and only if for any generator $x\in A$ in degree
$n+1$ the element $\xi(d_Ax)$ is a coboundary in $B$. (Note that
$d_Ax\in A_n$ so $\xi(d_Ax)$ is defined. We have
\[d_B\xi(d_Ax)=\xi(d_Ad_Ax)=0\]
which means that $\xi(d_Ax)$ is an $n+2$-cocycle in $B$. Since all
$n+1$-cocycles are coboundaries we see that $\xi$ can indeed be
extended to $A_{n+1}$. Using induction up the Postnikov tower of
$A$ we see that $\xi$ could be extended to a derivation $A\lra B$.

Injectivity: suppose that $\xi\in \Der^0(A,B)$ determines a
boundary in $\Der^0(A_n,B_n)$; we will then show that $\xi$ is a
boundary in $\Der^0(A,B)$. Indeed, considering $\xi$ as an element
in $\Der^0(A_n,B_n)$ we have
\[\xi=d_B\circ\eta+\eta\circ d_A\]
where $\eta$ is a derivation $A_n\lra B_n$ of degree $-1$. Take a
generator $x\in A$ in dimension $n+1$; we want to define $\eta(x)$
so that the following equality were true:
\[\xi(x)=d_B\circ\eta(x)+\eta\circ d_A(x).\]
For this, it is necessary and sufficient that
$(\xi-d_B\circ\eta)(x)$ be a coboundary in $B$. We have:
\begin{align*}d_A[(\xi-d_B\circ\eta)(x)]&=\xi(d_Ax)-d_A(\eta(d_Bx))\\&=\eta(d_B\circ
d_Bx)\\&=0
\end{align*}
In other words, $(\xi-d_B\circ\eta)(x)$ is an $n+1$-cocycle in
$B$. Since all $n+1$-cocycles are coboundaries we conclude that
$\xi$ restricted to $A_{n+1}$ is a coboundary in
$\Der^0(A_{n+1},B_{n+1})$. Induction up the Postnikov tower of $A$
finishes the proof.

\end{proof}
Let $X$ be a nilpotent space of finite type. Denote by $A$ its
minimal model. Then $A$ is an augmented commutative differential
graded algebra over $\mathbb{Q}$ which is a polynomial algebra on
$\pi_*(X)$. The differential on $A$ is a derivation $d_A:A\lra A$
of degree $1$ for which $d_A(A)=A_+\cdot A_+$. Here we denoted by
$A_+$ the set of elements in $A$ having positive degree.

Let $L:=\Der^* A$, the set of all graded derivations of $A$. The
differential on $L$ is given by the formula $d(\eta)=[\eta,d_A]$
where $\eta\in L$. The condition $d_A\circ d_A=0$ ensures that the
operator $[-,d_A]$ in $L$ has square $0$.

Since $A$ is a cofibrant object in the closed model category of
differential graded algebras the cohomology of $L$ represents the
derived functor of derivations. In other words, $H^*(L)\cong
H_{AQ}^*(A,A)$ as we saw in section 2.  The set of derivations of
degree $0$, that is $H_{AQ}^0(A,A)$ is then a conventional
(ungraded) Lie algebra.
\begin{rem}Consider the graded Lie algebra $B^*(L)$. The condition that
$A$ is minimal nilpotent ensures that every element in $B^*(L)$ is
nilpotent. In particular $B^0(L)$ is a nilpotent (ungraded) Lie
algebra over $\mathbb{Q}$.\end{rem}

We will now briefly recall the notion of homotopy in the category
of dga's restricting ourselves to self-maps. The details may be
found in \cite{BG}.

Let $A[t,dt]$ denote the differential graded algebra obtained from
$A$ by adjoining polynomial variables $t,dt$ subject to the
relation $(dt)^2=0$. The differential $d_{A[t,dt]}$ on $A[t,dt]$
is induced from the one on $A$. More precisely, denote by
$\partial_t$ the partial derivative with respect to $t$. Then for
$h\in A[t,dt]$ we have
\[d_{A[t,dt]}(h)=d_A(h)+(\partial_th)dt.\] There are two dga maps
$A[t,dt]\lra A$ given by $e_0:h\lra h|_{t=0}$ and $e_1:h\lra
h|_{t=1}$

 Let $F:A\lra A$ be a dga self map of $A$. Then
$F$ is said to be homotopic to the identity if there exists a dga
map $A\lra A[t,dt]$ such that its composition with $e_0$ is the
identity map on $A$ whereas its composition with $e_1$ is $F$.
Note that any map $A\lra A[t,dt]$ could be written as $F+Gdt$
where $F$ and $G$ are maps $A\lra A[t]$. The set of dga self-maps
of $A$ homotopic to the identity forms a normal subgroup in the
set of all automorphisms of $A$. The corresponding quotient group
is the group of \emph{homotopy self-equivalences} of $A$ and will
be denoted by $hAut(A)$. It is isomorphic to the group of homotopy
classes of homotopy self-equivalences of $X_{\mathbb{Q}}$, the
rationalization of the space $X$.
\begin{theorem}\label{harrison}
Let $X$ denote a nilpotent $CW$ complex which is either finite or
has a finite Postnikov tower, and let $A^*(X)$ denote its
Sullivan-de Rham model. The group $hAut(A^*(X))\cong
hAut(X_{\mathbb{Q}})$ is the group of $\mathbb{Q}$-points of an
affine algebraic group scheme over $\mathbb{Q}$ whose Lie algebra
is $H_{AQ}^0(A^*(X),A^*(X))$.
\end{theorem}
\begin{rem}Sullivan in \cite{Sul} sketched the proof of the fact that
$hAut(X_{\mathbb{Q}})$ is an algebraic group. The explicit
identification of its Lie algebra in terms of Andr\'e-Quillen
cohomology is new.
\end{rem}
\begin{proof}
First of all we replace $A^*(X)$ by its minimal model denoted by
$A$. Then $hAut(A)\cong hAut(A^*(X))$. By \ref{wellknown} and
\ref{wellknowntoo}, the case of finite CW complex reduces to that
of the finite Postnikov tower. Thus we may assume that $A$ is a
nilpotent minimal algebra with generators of bounded degree. Now
because of the finiteness assumptions the group of algebra
automorphisms of $A$ (i.e. not taking into account the
differential) is algebraic. Further the condition that an algebra
map $A\lra A$ commutes with $d_A$ is algebraic from which it
follows that the group $\mathcal{Z}$ of dga automorphisms of $A$
is also an algebraic group over $\mathbb{Q}$. It is obtained from
its \emph{reductive} part (coming from the quadratic part of the
differential $d_A$) by iterated extensions by the \emph{additive
group scheme}.

The Lie algebra of the group $\mathcal{Z}$ is just the set of
degree $0$ derivations of $A$ which commute with the differential
in $A$. Therefore this is the Lie algebra $Z_0(L)$ of zero degree
cocycles of $L$. The normal Lie subalgebra $B^0(L)$ of $Z^0(L)$ is
nilpotent and there is a corresponding normal subgroup
$\mathcal{B}$ in the algebra self-maps of $A$ obtained by
exponentiating $B^0(L)$. We will show that this subgroup consists
precisely of those self-maps which are homotopic to the identity.
Thus, $hAut(A)$ is the quotient of an affine algebraic group by a
normal affine algebraic subgroup and is thus affine algebraic.

Granting this for a moment, we see that the Lie algebra of
$hAut(A)=\mathcal{Z/B}$ is the quotient Lie algebra
$Z_0(L)/B_0(L)$ which is precisely $H^0_{AQ}(A,A)\cong
H^0_{AQ}(A^*(X),A^*(X))$ by the results of Section 2.

Let $F+Gdt:A\lra A[t,dt]$ be a homotopy for which $F|_{t=0}=id_A$.
First examine the condition that it is a dga map. For $h_1,h_2\in
A$ we have:
\begin{align*}(F+Gdt)(h_1h_2)=&F(h_1h_2)+G(h_1h_2)dt\\
=&F(h_1)F(h_2)+[F(h_1)G(h_2)+G(h_1)F(h_2)]dt\\=&(F+Gdt)(h_1)(F+Gdt)(h_2).
\end{align*}
So we get two conditions:\begin{enumerate}\item
$F(h_1h_2)=F(h_1)F(h_2)$ which means simply that $F$ is an algebra
map and \item $G(h_1h_2)=F(h_1)G(h_2)+G(h_1)F(h_2)$
\end{enumerate}
The second condition means that $G$ is an $F$-derivation. Setting
$t=0$ in the second equation we get
\[G|_{t=0}(h_1h_2)=h_1G(g_2)+G(h_1)h_2.\]
In other words $G|_{t=0}$ is a usual derivation of $A$.

Next, the condition that $F+Gdt$ is a map of \emph{differential}
algebras means that
\[(F+Gdt)\circ d_A=d_{A[t,dt]}\circ(F+Gdt).\]
Applying both sides of this equation to $h\in A$ and equating
coefficients at $t$ and $dt$ we get two
identities:\begin{enumerate}\item $F\circ d_A=d_A\circ F$\item
$\partial_tF=G\circ d_A+d_A\circ G=[G,d_A]$.
\end{enumerate}
The first condition above simply means that $F=F(t)$ commutes with
the differential; i.e. determines a family of maps of complexes
$A\lra A$. Setting $t=0$ in the second equation we get
$\partial_tF|_{t=0}=[G|_{t=0},d_A]$ which means that
$\partial_tF|_{t=0}$ is a coboundary in $L$.

Now let $F_1$ be an element in $\mathcal{B}$. Then $F$ could be
represented as $F_1=\exp([G_0,d_A])$ where $G_0$ is a derivation
of $A$. Let $F=F(t)=\exp(t[G_0,d_A])$ and
$G=G(t)=G_0\exp(t[G_0,d_A])$. Then $(F(t),G(t))$ is the homotopy
from $F_1$ to $id$. Indeed, $F(t)$ is an algebra map for each $t$
and $G(t)$ is an $F(t)$-derivation.

Conversely, suppose that $(F(t),G(t))$ is a homotopy from $id $ to
$F(1)$, we need to show that $F(1)$ belongs to $\mathcal{B}$. Here
we will use the fact that $A$ is nilpotent and minimal but our
proof extends under the condition that $F_1$ is homotopic to the
identity through automorphisms, that is $F(t)$ is an automorphism
for all $t$. This condition follows automatically for minimal
algebras since a weak equivalence between minimal nilpotent
algebras is necessarily an isomorphism.

We can write \[F(t)=Id_A+\sum_{j=1}^\infty F_jt^j.\] Since
$F:A\lra A[t,dt]$ we must have that the sum is locally finite in
the sense that for $a\in A$, $F_j(a)=0$ for $j\gg 0$. $F^{-1}$
exists formally: \[F^{-1}(t)=Id_A+\sum_{i=1}^\infty
(-1)^i(\sum_{j=1}^\infty F_jt^j)^i\] and this infinite sum is
locally finite by the condition above.

Since $\partial_tF=[G,d_A]$ we have, taking into account that $F$
commutes with $d_A$:
\begin{align*}(\partial_tF)F^{-1}=&[G,d_A]F^{-1}\\=&G\circ d_A \circ
F^{-1}-d_A\circ G\circ F^{-1}\\=&[GF^{-1},d_A].\end{align*} Since
$G(t)$ is an $F(t)$-derivation $GF^{-1}(t)$ will be a usual
derivation. Noting that \[(\partial_tF(t))F^{-1}(t)=\partial_t
\log F(t)\] we have the following equation: \[ \partial_t\log
F(t)=[GF^{-1},d_A].\]
Therefore \begin{align*}F(1)&=\exp(\int_0^1[GF^{-1},d_A]dt)\\
&=\exp[\int_0^1GF^{-1}dt,d_A]
\end{align*}
Here the integral is carried out formally. Moreover,
$\int_0^1GF^{-1}dt$ is a locally finite expression and still a
derivation. Therefore $F(1)$ is contained in $\mathcal{B}$.

\end{proof}
\begin{rem}One naturally wonders whether there is a relationship
between the whole complex $C^*_{AQ}(A^*(X),A^*(X))$ together with
its Gerstenhaber bracket and the space $Aut(X)$ of homotopy
self-equivalences of $X$. Schlessinger and Stasheff, \cite{SS1}
provide an affirmative answer to this question in the case when
$X$ is simply-connected. Namely, they show that the complex
consisting of derivations of $A$ lowering the degree by $k>1$
serves as a Lie model for the universal covering of the space
$BAut(X)$. That implies that the Whitehead product in
$\pi_{>1}BAut(X)$ and the Gerstenhaber bracket in
$H^{<0}_{AQ}(A^*(X),A^*(X))$ agree whereas our theorem compares
$\pi_{1}BAut(X)$ and $H^{0}_{AQ}(A,A)$. It seems likely that the
theorem of Schlessinger-Stasheff could be extended to the
nilpotent case as well.
\end{rem}
\begin{rem}\label{Lie}
One can also express the tangent Lie algebra to $hAut(X_{Q})$ in
terms of the Quillen or dg Lie algebra model of $X$, at least when
$X$ is simply-connected (cf. \cite{HFT} concerning Quillen
models.) Namely, if $L(X)$ is a Quillen model for $X$ consider the
Chevalley-Eilenberg complex $C^*_{CE}(L(X),\mathbb{Q})$ computing
the Lie algebra cohomology of $L$ with trivial coefficients. Then
$C^*_{CE}(L(X),\mathbb{Q})$ is a cofibrant dga and serves as a
Sullivan model for $X$. Furthermore the dg Lie algebra of
derivations of $C^*_{CE}(L(X),\mathbb{Q})$ coincides up to a shift
of dimensions with the complex $C^*_{CE}(L(X),L(X))$ computing the
Lie algebra cohomology of $L(X)$ with coefficients in itself. The
complex $C^{*-1}_{CE}(L(X),L(X))$ carries a dg Lie algebra
structure so that $H^1_{CE}(L(X),L(X))$, its cohomology in the
total degree $1$ is an ungraded Lie algebra. We conclude that the
Lie algebra of $hAut(X_{Q})$ is isomorphic to
$H^1_{CE}(L(X),L(X))$.
\end{rem}
Our next result is concerned with the more general problem of
computing $\pi_iF(X,Y)$ for $i\geq 1$. Here $X,Y$ are rational
nilpotent $CW$-complexes of finite type and $F(X,Y)$ is the space
of continuous maps from $X$ into $Y$. Recall that $A^*(X)$ and
$A^*(Y)$ denote the Sullivan-de Rham algebras of $X$ and $Y$
respectively. In view of the previous remark it is natural to
expect that the homotopy groups of $F(X,Y)$ are expressed in terms
of the Andr\'e-Qullen cohomology of $A^*(Y)$ with values in
$A^*(X)$. When $X=Y$ and $X$ is simply-connected the corresponding
statement was proved by Schlessinger and Stasheff by combining
their deformation theory with the fact that $BAut(X)$ is the
classifying space for $X$-fibrations. This machinery is not
available for $X\neq Y$ and we use Proposition
\ref{characterization} instead. Note also that the Gerstenhaber
bracket no longer exists on the complex $C^*_{AQ}(A,B)$ for $A\neq
B$.

Assume that there is given a map $f:X\lra Y$ which determines a
basepoint in $F(X,Y)$.  The map $f$ induces a map $A^*(Y)\lra
A^*(X)$ making $A^*(X)$ into a differential graded module over
$A^*(Y)$.
\begin{theorem} There is an isomorphism of sets
\[\pi_n(F(X,Y),f)\cong H^{-n}_{AQ}(A^*(Y),A^*(X)).\]
If $n\geq 2$ then this is an isomorphism of abelian groups.
\end{theorem}
\begin{proof}Consider the de Rham algebra  $A^*(X\times S^n)$ where
$n\geq 1$. Clearly $A^*(X\times S^n)$ is weakly equivalent as a
commutative dga to the dga $A^*(X)\ltimes A^*(X)[n]$ where
$A^*(X)[n]$ is the square-zero ideal which is isomorphic as a
$\mathbb{Q}$-vector space to $A^*(X)$ and whose grading is shifted
by $n$.

Denoting by $F_*(?,?)$ the function space between \emph{pointed}
topological spaces we have an obvious homotopy fibre sequence:
\[F_*(S^n,F(X,Y))\lra F(S^n,F(X,Y))\lra F(*,F(X,Y))=F(X,Y)\]
where the second arrow is induce by the inclusion of the basepoint
in $S^n$. (The basepoint of $F(X,Y)$ is $f$ and we take
$F_*(S^n,F(X,Y))$ to be the fibre \emph{over $f$}.) Standard
adjunction gives the homotopy fibre sequence
\begin{equation}\label{fibre}
F_*(S^n,F(X,Y))\lra F(X\times S^n,Y)\lra F(X,Y).\end{equation}
Consider now the category ${\mathcal{T}_X}$ of spaces \emph{under}
$X$, i.e. spaces supplied with a map from $X$ with morphisms being
the obvious commutative triangles. This is a topological closed
model category and $X\times S^n,Y\in {\mathcal T}_X$. Note that
$X\times S^n$ is a cogroup object in the homotopy category of
${\mathcal T}_X$, abelian when $n>1$. Likewise $A^*(X\times S^n)$
is (equivalent to) an abelian group object in
$\mathcal{C}_{A^*(X)}$. For $n>1$ the cogroup structure on
$X\times S^n$ corresponds to the group structure on $A^*(X\times
S^n)$ under the Sullivan-de Rham functor
$A^*:\mathcal{T}_X\lra\mathcal{C}_{A^*(X)}$.

It follows that  the space $F_*(S^n,F(X,Y))$ is weakly equivalent
to the function space in $\mathcal{T}_X$ from $X\times S^n$ to
$Y$. Therefore, denoting by $[-,-]_{\mathcal{T}_X}$ the set of
homotopy classes of maps in $\mathcal{T}_X$ we have:
\[\pi_nF_*(X,Y)\cong [X\times S^n,Y]_{\mathcal{T}_X}.\]
Since the homotopy category of finite type rational spaces in
$\mathcal{T}_X$ is anti-equivalent to the homotopy category of
finite type dga's \emph{over} $A^*(X)$ we have

\begin{align*}\pi_nF_*(X,Y)&\cong [A^*(Y),A^*(X\times
S^n)]_{\mathcal{C}_{A^*(X)}}
\\ &\cong
[A^*(Y),A^*(X)\ltimes A^*(X)[n]]_{\mathcal{C}_{A^*(X)}}
\\
&\cong H_{AQ}^0(A^*(Y,)A^*(X)[n])~\mbox{(by Proposition
\ref{characterization})}\\&\cong H_{AQ}^{-n}(A^*(Y),A^*(X))
\end{align*}
\end{proof}
\begin{rem}If $X$ is a point then the space of maps $X\lra Y$ is
simply $Y$. The previous theorem then gives an identification of
$\pi_{*}(Y)$ with $H_{AQ}^{*}(A^*(Y),\mathbb{Q})$. This is not
surprising since by Theorem \ref{derived}
$H_{AQ}^{*}(A^*(Y),\mathbb{Q})$ is identified with the homology of
the complex $\Der^*(M(A^*(Y)))$.  Here $M(A^*(Y))$ denotes the
minimal model of $Y$. Clearly this homology is simply the dual
space of indecomposables in $M(A^*(Y))$ so we recover a standard
result in rational homotopy theory.

More generally, suppose that the map $f:X\lra Y$ is homotopic to
the trivial map. Thus, $A^*(X)$ is an $A^*(Y)$-module via the
augmentation map $A^*(Y)\lra \mathbb{Q}$. Then it is easy to see
that there is an isomorphism of complexes
\[C^*_{AQ}(A^*(Y),A^*(X))\cong
C^*_{AQ}(A^*(Y),\mathbb{Q})\hat{\otimes} A^*(X).\] Here
$\hat{\otimes}$ denotes completed tensor product. Therefore
\begin{align*}H^*_{AQ}(A^*(Y),A^*(X))\cong&H^*_{AQ}(A^*(Y),\mathbb{Q})\hat{\otimes}
H^*(X)\\\cong&\pi_*(Y)\hat{\otimes} H^*(X).\end{align*} So we
recover an isomorphism obtained in \cite{BS}:
\[\pi_n(F(X,Y),f)\cong \prod_{k=1}^\infty\pi_k(Y)\otimes
H^{k-n}(X).\]
\end{rem}
\begin{rem}In the case $X$ and $Y$ are simply-connected the
considerations similar to those in Remark \ref{Lie} yield a
bijection of sets for $n=1$ and an isomorphism of abelian groups
for $n>1$:
\[\pi_nF(X,Y)\cong H^{1-n}(L(X),L(Y))\]
where $L(X)$ and $L(Y)$ are Quillen models of $X$ and $Y$
respectively. It seems likely that the last isomorphism continues
to hold under no finiteness assumptions on the rational spaces $X$
and $Y$.\end{rem}

\end{document}